\documentclass[12TP,draft]{article}

\begin{document}

\begin{center}
\LARGE\noindent\textbf{ Cycles of each even lengths in balanced bipartite digraphs }\\

\end{center}
\begin{center}
\noindent\textbf{Samvel Kh. Darbinyan }\\

Institute for Informatics and Automation Problems, Armenian National Academy of Sciences

E-mails: samdarbin@ipia.sci.am\\
\end{center}

\textbf{Abstract}

Let $D$ be a strongly connected balanced bipartite directed graph of order $2a\geq 4$. Let $x,y$ be distinct vertices in $D$. $\{x,y\}$ dominates a vertex $z$ if $x\rightarrow z$ and $y\rightarrow z$; in this case, we call the pair $\{x,y\}$ dominating.  In this paper we prove: 

(i). {\it If  $a\geq 4$ and $ max\{d(x), d(y)\}\geq 2a-1$ for every dominating pair of vertices $\{x,y\}$, $D$ then contains a cycle of length $2a-2$ or $D$ is a directed cycle.}

(ii). {\it If  $D$ contains a cycle of length $2a-2\geq 6$ and $max \{d(x), d(y)\}\geq 2a-2$ for every dominating pair of vertices $\{x,y\}$, then for any $k$, $1\leq k\leq a-1$, $D$ contains a cycle of length $2k$.}

(iii). {\it If $a\geq 4$ and $ max\{d(x), d(y)\}\geq 2a-1$ for every dominating pair of vertices $\{x,y\}$,  then for every $k$, $1\leq k\leq a$,   $D$ contains a cycle of length $2k$ unless $D$ is isomorphic to only one exceptional digraph of order eight.}\\


\section {Introduction} 

We consider directed graphs (digraphs) in the sense of \cite{[3]}. For convenience of the reader terminology and notations will be given in details in section 2. A cycle (path) is called Hamiltonian  if it  includes all the vertices of $D$. A digraph $D$ is Hamiltonian  if it contains a  Hamiltonian cycle and is pancyclic if it contains a cycle of length $k$ for any $3\leq k\leq n$,
 where $n$ is the order of $D$. A balanced bipartite digraph of order $2m$ is even pancyclic  if it contains a cycle of length $2k$ for any $k$, $2\leq k\leq m$. 

 Bondy suggested  \cite{[6]} the following interesting  "metaconjecture":\\

 {\it Almost any non-trivial condition  of a graph (digraph) which that the graph (digraph) is Hamiltonian  also implies that the graph (digraph) is pancyclic. (There may be  a "simple" family of exceptional graphs (digraphs)).}\\

 There are numerous sufficient conditions for the existence of a Hamiltonian cycle in a graph (digraph) for which "metaconjecture" of Bondy is true.

In this note, we will concerted with the degree conditions. 
Let us recall the following well-known degree conditions (theorems 1.1-1.4) that guarantee that a digraph is Hamiltonian.
\\
 \textbf{Theorem 1.1} (Ghouila-Houri \cite{[12]}). {\it Let $D$ be a strongly connected digraph of order $n\geq 3$. If $d(x)\geq n$ for all vertices $x\in V(D)$, then $D$ is Hamiltonian.}\\

 \textbf{Theorem 1.2} (Woodall \cite{[20]}). {\it Let $D$ be a digraph of order $n\geq 3$. If $d^+(x)+d^-(y)\geq n$ for all pairs of vertices $x$ and $y$ such that there is no arc from $x$ to $y$, then $D$ is Hamiltonian.}\\
 
\textbf{Theorem 1.3} (Meyniel \cite{[16]}). {\it Let $D$ be a strongly connected digraph of order $n\geq 2$. If $d(x)+d(y)\geq 2n-1$ for all pairs of non-adjacent vertices in $D$, then $D$ is Hamiltonian.}\\

Notice that Meyniel's theorem is a common generalization of Ghouila-Houri's and Woodall's theorems. For a short proof of Theorem 1.3, see \cite{[7]}.

C. Thomassen \cite{[18]} (for $n=2k+1$) and S. Darbinyan \cite{[9]}
(for $n=2k$) proved the following theorem below. 
Before stating it (for any integer $m\geq 2$)  we need to introduce some additional notations.\\

Let $H(m,m)$ denote the set of digraphs $D$ of order $2m$ with vertex set $A\cup  B$ such that $\langle A\rangle\equiv \langle B\rangle \equiv K^*_m$, there is no arc from $B$ to $A$, $d^+(x,B)\geq 1$ and $d^-(y,A)\geq 1$ for every vertices $x\in A$ and $y\in B$.

Let $H(m,m-1,1)$ denote the set of digraphs $D$ of order $2m$ with vertex set $A\cup  B\cup \{a\}$ such that $|A|=|B|+1=m$, $\langle B\cup \{a\}\rangle \subseteq K^*_m$ (i.e., $\langle B\cup \{a\}\rangle$  is an arbitrary digraph) the subdigraph $\langle A\rangle $ has no arc, $D$ contains all possible arcs between $A$ and $B$  and either $N^-(a)=B$ and $A\subseteq N^+(a)$, or $N^+(a)=B$ and $A\subseteq N^-(a)$.

Let $H(2m)$ denote a digraph of order $2m$ with vertex set $A\cup B \cup \{x,y\}$ such that $\langle A\cup \{x\}\rangle \equiv\langle B\cup \{y\}\rangle \equiv K^*_m$, there is no are between $A$ and $B$, $H(2m)$ also contains all arcs of the form $ya$, $bx$, where $a\in A$ and $b\in B$, and either the arc $xy$ or both arcs $xy$ and $yx$. 

Let $D_6$ be a digraph with vertex set $\{x_1,x_2,\ldots ,x_5,x\}$  and arc set 
$$
\{x_ix_{i+1} \,/ 1\leq i\leq 4 \}\cup \{xx_i/\,1\leq i\leq 3\} \cup \{x_1x_5,x_2x_5,x_5x_1,x_5x_4,x_3x_2,x_3x,x_4x_1,x_4x \}.$$ 
By $D'_6$ we denote a digraph obtained from $D_6$ by adding the arc  $x_2x_4$. 
 Note that the digraphs  $D_6$ and $D'_6$ both are not Hamiltonian and each  of $D_6$ and $D'_6$ contains a cycle of length 5. \\ 

\noindent\textbf{Theorem 1.4} ( Thomassen \cite{[18]}, Darbinyan \cite{[9]}). {\it If $D$ is a digraph of order $n\geq 5$ with minimum degree at least $n-1$ and with minimum semi-degree at least $n/2-1$. Then $D$ is Hamiltonian unless 

(i) $D$ is isomorphic to $D_5$ or $D_7$ or
$[(K_{m}\cup K_{m})+K_1]^*$ or $K^*_{m,m+1}\subseteq D \subseteq [K_{m}+\overline K_{m+1}]^*,$ if $n=2m+1$;

(ii) $D\in H(m,m)\cup H(m,m-1,1)\cup  \{ H(2m),\,H'(2m),\, D_6,\, D'_6,\overleftarrow{D_6},\, \overleftarrow {D'_6}\}$, if $n=2m$. 
(The digraphs $D_5$ and $D_7$ are well-known and for  definitions of their,  see, for example, \cite{[18]}).} \fbox \\\\

 In \cite{[15], [17], [8], [10]} it was shown that if a digraph $D$ satisfies the condition one of Theorems 1.1-1.4, respectively, then $D$ also is pancyclic (unless some extremal cases which are characterized). 

 Characterizations of even pancyclic bipartite tournaments was given in \cite{[5]} and \cite{[21]}. A characterizations of pancyclic ordinary $k$-partite ($k\geq 3$) tournaments was established in \cite{[13]}. A characterizations of pancyclic ordinary complete $k$-partite ($k\geq 3$) digraphs was derived in \cite{[14]}. Amar and Manoussakis \cite{[2]} gave several sufficient conditions on the half-degrees of a bipartite digraph for the existence of cycles and paths of various lengths.

Each of  Theorems 1.1-1.4  imposes a degree condition on all pairs of nonadjacent vertices (or on all vertices). In   \cite{[2], [4]} were shown some sufficient conditions for hamiltonicity of digraphs in which the degree conditions require only for some pairs of nonadjacent vertices. Let us recall of them only the following theorem. \\

\textbf{Theorem 1.5} (Bang-Jensen, Gutin, H.Li \cite{[4]}). {\it Let $D$ be a strongly connected digraph of order $n\geq 2$. Suppose that $min\{d(x),d(y)\}\geq n-1$ and  $d(x)+d(y)\geq 2n-1$ for any pair of nonadjacent vertices $x,y$ with a common in-neighbor. Then $D$ is Hamiltonian.}\\

An analogue of Theorem 1.5  for bipartite digraphs was given by R. Wang \cite{[19]} and recently strengthened by the author \cite{[11]}.

\textbf{Theorem 1.6} (R. Wang \cite{[19]}). {\it Let $D$ be a strongly connected balanced bipartite digraph of order $2a$, where $a\geq 1$. Suppose that, for every dominating pair of vertices $\{x,y\}$, either $d(x)\geq 2a-1$ and $d(y)\geq a+1$ or $d(y)\geq 2a-1$ and $d(x)\geq a+1$. Then $D$ is Hamiltonian}.\\

Before stating the next theorem we need to define a digraph of order eight.\\

\textbf{Definition.} Let $D(8)$ be a  bipartite digraph  with partite sets $X=\{x_0,x_1,x_2,x_3\}$ and 
$Y=\{y_0,y_1,y_2,y_3\}$, and let $D(8)$ contains the arcs $y_0x_1$, $y_1x_0$, $x_2y_3$, $x_3y_2$ and all the arcs of the following 2-cycles: 
$x_i\leftrightarrow y_i$, $i\in [0,3]$, $y_0\leftrightarrow x_2$, $y_0\leftrightarrow x_3$, $y_1\leftrightarrow x_2$ and it $y_1\leftrightarrow x_3$, and it contains no other arcs. 

It is not difficult to check that $D(8)$ is strongly connected, $max \{d(x), d(y)\}\geq 2a-2$ for every  pair of vertices $\{x,y\}$ with a common in-neighbor and is not Hamiltonian.\\

We also need the following definition.

\textbf{Definition.} Let $D$ be a balanced bipartite digraph of order $2a$, where $a\geq 2$, and let  $k$ be an integer.
 We say that
 $D$ satisfies condition $B_k$ if  for every dominating pair of vertices  $x$ and $y$, $max\{d(x),d(y)\}\geq 2a-2+k$ .\\

\textbf{Theorem 1.7} (Darbinyan \cite{[11]}). {\it Let $D$ be a strongly connected balanced bipartite digraph of order $2a$, where $a\geq 4$. Suppose that  $max\{d(x),d(y)\}\geq 2a-1$, for every dominating pair of vertices $\{x,y\}$. Then either $D$ is Hamiltonian or isomorphic to the digraph $D(8)$}.\\

Motivated by the Bondy's "metaconjecture", it is natural to set the following problem: 

{\it Characterize those digraphs which satisfy the conditions of Theorem 1.6 (1.7) but are not even pancyclic.}

In this note we prove the following theorems. 

 \textbf{Theorem 1.8.} {\it Let $D$ be a strongly connected balanced bipartite digraph of order $2a\geq 8$ 
with partite sets $X$ and $Y$. Suppose that $D$ satisfies condition $B_1$, i.e., $ max\{d(x), d(y)\}\geq 2a-1$ for every dominating pair of vertices $\{x,y\}$.
  Then  $D$ contains a cycle of length $2a-2$ or $D$ is a directed cycle.}\\

\textbf{Theorem 1.9.} {\it Let $D$ be a strongly connected balanced bipartite digraph of order $2a\geq 8$ with partite sets $X$ and $Y$. If  $D$ contains a cycle of length $2a-2$ and satisfies condition $B_0$, i.e., $max \{d(x), d(y)\}\geq 2a-2$ for every dominating pair of vertices $\{x,y\}$,
 then  $D$ contains a cycle of length $2k$ for all $k$, where $1\leq k\leq a-1$.}\\

\textbf{Theorem 1.10.} {\it Let $D$ be a strongly connected balanced bipartite digraph of order $2a\geq 8$ with partite sets $X$ and $Y$.  If $D$ is not a directed cycle and $max\{d(x), d(y)\}\geq 2a-1$ for every dominating pair of vertices $\{x,y\}$,  then either $D$ contains  cycles of all even lengths less than  equal to $2a$ or $D$ is isomorphic to the digraph $D(8)$.}\\

Of course, Theorem 1.10 is an immediate corollary of Theorems 1.7-1.9.

\section {Terminology and Notations}

  In this paper we consider finite digraphs without loops and multiple arcs. Terminology  and notations not described below follow \cite{[2]}. 
For a digraph $D$, we denote
  by $V(D)$ the vertex set of $D$ and by  $A(D)$ the set of arcs in $D$. The order of $D$ is the number
  of its vertices. 
The arc of a digraph $D$ directed from
   $x$ to $y$ is denoted by $xy$ or $x\rightarrow y$ (we also say that $x$ dominates $y$ or $y$ is an out-neighbour of $x$ and $x$ is an in-neighbour of $y$), and $x\leftrightarrow y$ denotes that $x\rightarrow y$ and $y\rightarrow x$ ($x\leftrightarrow y$ is called  2-cycle). If $x\rightarrow y$ and $y\rightarrow z$ we write $x\rightarrow y\rightarrow z$.
If there is no arc from $x$ to $y$ we shall use the notation $xy\notin A(D)$.
For disjoint subsets $A$ and  $B$ of $V(D)$  we define $A(A\rightarrow B)$ \,
   as the set $\{xy\in A(D) / x\in A, y\in B\}$ and $A(A,B)=A(A\rightarrow B)\cup A(B\rightarrow A)$. If $x\in V(D)$
   and $A=\{x\}$ we sometimes write $x$ instead of $\{x\}$. If $A$ and $B$ are two disjoint subsets of $V(D)$ such that every
   vertex of $A$ dominates every vertex of $B$, then we say that $A$ dominates $B$, denoted by $A\rightarrow B$. $A\mapsto B$ means that $A\rightarrow B$ and there is no arc from $B$ to $A$. Similarly, $A\leftrightarrow B$ means that 
$A\rightarrow B$ and $B\rightarrow A$.
The out-neighborhood of a vertex $x$ is the set $N^+(x)=\{y\in V(D) / xy\in A(D)\}$ and $N^-(x)=\{y\in V(D) / yx\in A(D)\}$ is the in-neighborhood of $x$. Similarly, if $A\subseteq V(D)$, then $N^+(x,A)=\{y\in A / xy\in A(D)\}$ and $N^-(x,A)=\{y\in A / yx\in A(D)\}$. 
The out-degree of $x$ is $d^+(x)=|N^+(x)|$ and $d^-(x)=|N^-(x)|$ is the in-degree of $x$. Similarly, $d^+(x,A)=|N^+(x,A)|$ and $d^-(x,A)=|N^-(x,A)|$. The degree of the vertex $x$ in $D$ is defined as $d(x)=d^+(x)+d^-(x)$ (similarly, $d(x,A)=d^+(x,A)+d^-(x,A)$). The subdigraph of $D$ induced by a subset $A$ of $V(D)$ is denoted by $D\langle A\rangle$ or $\langle A\rangle$ brevity. The path (respectively, the cycle) consisting of the distinct vertices $x_1,x_2,\ldots ,x_m$ ( $m\geq 2 $) and the arcs $x_ix_{i+1}$, $i\in [1,m-1]$  (respectively, $x_ix_{i+1}$, $i\in [1,m-1]$, and $x_mx_1$), is denoted by  $x_1x_2\cdots x_m$ (respectively, $x_1x_2\cdots x_mx_1$). The length of a cycle or a path is the number of its arcs.
We say that $x_1x_2\cdots x_m$ is a path from $x_1$ to $x_m$ or is an $(x_1,x_m)$-path. The length of a cycle or a path is the number of its arcs.
 
If $P$ is a path containing a subpath from $x$ to $y$ we let $P[x,y]$ denote that subpath.
 Similarly, if $C$ is a cycle containing vertices $x$ and $y$, $C[x,y]$ denotes the subpath of $C$ from $x$ to $y$.
 Given a vertex $x$ of a path $P$ or a cycle $C$, we denote by $x^+$ (respectively, by $x^-$) the successor (respectively, the predecessor) of $x$ (on $P$ or $C$), and in case of ambiguity, we precise $P$ or $C$ as a subscript (that is $x^+_P$ ...). 

A digraph $D$ is strongly connected (or, just, strong) if there exists a path from $x$ to $y$ and a path from $y$ to $x$ for every pair of distinct vertices $x,y$.
    
Two distinct vertices $x$ and $y$ are adjacent if $xy\in A(D)$ or $yx\in A(D) $ (or both). For integers $a$ and $b$, $a\leq b$, let $[a,b]$  denote the set of
all the integers which are not less than $a$ and are not greater than
$b$. Let $C$ be a non-Hamiltonian cycle in digraph $D$. An $(x,y)$-path $P$ is a $C$-bypass if $|V(P)|\geq 3$, $x\not=y$ and $V(P)\cap V(C)=\{x,y\}$. The length of the path $C[x,y]$ is the gap of $P$ with respect to $C$.

A digraph $D$ is a bipartite  if there exists a partition $X$, $Y$ of $V(D)$ into two partite sets such that every arc of $D$ has its end-vertices in different partite sets. 
It is called balanced if $|X|=|Y|$. The underlying graph of a digraph $D$ is denoted by $UG(D)$, it contains an edge  $xy$ if $x\rightarrow y$ or $y\rightarrow x$ (or both).

\section { Preliminaries }

Let us recall some results (Lemmas 3.1-3
.3) which will be used in this paper.\\

  \textbf{Lemma 3.1} (\cite{[11]}). {\it Let $D$ be a strongly connected balanced bipartite digraph of order $2a\geq 8$. If $D$ satisfies condition $B_1$, then the following holds:

(i) $UG(D)$ is 2-connected;

(ii) if $C$ is a cycle of length $m$, $2\leq m\leq 2a-2$, then $D$ contains a $C$-bypass.}\\

\textbf{Lemma 3.2} (\cite{[11]}). {\it Let $D$ be a strongly connected balanced bipartite digraph of order $2a\geq 8$ other than the directed cycle of length $2a$. If $D$ satisfies condition $B_0$, then $D$ has a non-Hamiltonian cycle of length at least four. }\\

\textbf{Lemma 3.3} (\cite{[1]}). {\it Let $D$ be a bipartite digraph of order $n$ which contains a cycle $C$ 
of length $2b$, where $2\leq 2b\leq n-1$. Let $x$ be a vertex not contained in $C$. If $d(x,V(C))\geq b+1$, then  $D$ contains cycles of every even length $m$, $2\leq m\leq 2b$, through $x$.}\\

Now we prove the following lemma.\\

\textbf{Lemma 3.4.} {\it Let $D$ be a strongly connected balanced bipartite digraph of order $2a\geq 8$  with partite sets $X$ and $Y$. Assume that $D$  satisfies condition $B_0$. Let $C=x_1y_1x_2y_2\ldots x_ky_kx_1$ be a longest non-Hamiltonian cycle in $D$, where $k\geq 2$, $x_i\in X$ and $y_i\in Y$. If in $D$ there exists a $C$-bypass whose gap with respect to $C$ is equal to one, then $k=a-1$, i.e., the longest non-Hamiltonian cycle in $D$ has length equal to $2a-2$. }\\

\textbf{Proof of Lemma 3.4}. Without loss of generality, we assume that $P:=x_1u_1u_2\ldots u_sy_1$ ($s\geq 1$) is a $C$-bypass. 
Suppose that the lemma is not true, that is $k\leq a-2$. Then $u_1\in Y$ and $u_s\in X$ and $R:=V(D)\setminus V(C)=\{u_1,u_2, \ldots , u_s\}$. Note that $|R|=s\geq 4$. Since $C$ is a longest non-Hamiltonian cycle in $D$, it is not difficult to show that  
$$
d^+(u_1,\{u_3,u_4,\ldots, , u_s\})=d^-(u_s,\{u_1,u_2,\ldots , u_{s-2}\})=d^+(x_1,\{u_2,u_3,\ldots , u_s\})=0, \eqno (1)
$$
and the following arcs
$$
u_{s-1} x_2,\, y_{k} u_2,\, x_{k} u_1,\, u_{s} y_2,\, u_{3} x_2, x_{1} u_3,\, y_{k-1} x_1, \, y_{1} x_3, \,  u_{2}y_1 \eqno (2)
$$
 are not in $A(D)$ (for otherwise, $D$ contains a non-Hamiltonian cycle longer than $C$). Note that $\{u_s,x_1\}$ is a dominating pair. Therefore, by condition $B_0$ we have 
$$
 max\{d(x_1),d(u_s)\}\geq 2a-2. \eqno (3)
$$

If $|R|=s\geq 6$, then from (1) and (2) we have $d^-(u_s,\{u_1,u_3\})=d^+(x_1, \{u_3,u_5\})=0$ and the arcs
  $u_sy_2$ and 
 $y_{k-1} x_1$ are not in $A(D)$. Therefore, $d(u_s)$ and $d(x_1)\leq 2a-3$, which contradicts (3).

Assume therefore that $|R|=4$ (i.e., $s=4$) and consider the following two cases.

\textbf{Case 1}. {\it The vertices $u_1$ and $u_4$ are not adjacent.}

Then
$$
d(u_1)\leq 2a-3 \quad \hbox{and} \quad d(u_4)\leq 2a-3   \eqno (4)
$$
since $x_{k} u_1$ and  $u_{4} y_2 \notin A(D)$ by (2).  From (3) and the second inequality of (4) it follows that $d(x_1)\geq 2a-2$. 
By (2) we have  $x_{1} u_3$ and  $y_{k-1} x_1\notin A(D)$. This together with $d(x_1)\geq 2a-2$ implies that $\{u_1,u_3,y_1\}\rightarrow x_1$. 
In particular, $\{u_1,u_3\}$ is a dominating pair. From the first inequality of (4) and condition $B_0$ it follows that $d(u_3)\geq 2a-2$. On the other hand, by (2) we have
 $x_{1}u_3$ and $u_{3}x_2\notin A(D)$. Hence, $u_{4}\rightarrow u_3\rightarrow u_2$, $x_k\rightarrow u_3$. 
Thus, $\{u_2,u_4\}\rightarrow u_3$, which means that $\{u_2,u_4\}$ is a dominating pair.
 Therefore, by   
 condition $B_0$, $d(u_2)\geq 2a-2$ since $d(u_4)\geq 2a-3$ by (4). From $d(u_2)\geq 2a-2$ and  $y_{k} u_2\notin A(D)$, $u_{2}y_1\notin A(D)$ (by (2)) it follows that $u_2\rightarrow u_1$. Thus, $x_ku_3u_2u_1x_1y_1\ldots y_{k-1}x_k$ is a cycle of length $2a-2$, which is a contradiction.

\textbf{Case 2}. {\it The vertices $u_1$ and $u_4$ are  adjacent.}

Then, by (1), we have $u_4\mapsto u_1$. We divide this case into two  subcases.

\textbf{Subcase 2.1}. $d(x_1)\leq 2a-3$.

Then $d(u_4)\geq 2a-2$ because of (3) (recall that $s=4$). Then, since, by (1) and (2), the arcs $u_1u_4$ and $ u_4 y_2$ are not in $A(D)$,  it follows  that
$$
d(u_4)= 2a-2, \quad y_1\rightarrow u_4\rightarrow u_3 \quad \hbox{and} \quad y_2\rightarrow u_4. \eqno (5)
$$
If $u_3\rightarrow x_3$, then by (5) we have that  $y_1u_4u_1u_2u_3x_3\ldots y_1$ is a cycle of length $2a-2$, a contradiction. 
Assume therefore that $u_3 x_3\notin A(D)$ (possibly $x_3=x_1$). This together with $x_{1} u_3\notin A(D)$ and 
$u_{3}x_2\notin A(D)$ (by (2)) implies that $d(u_3)\leq 2a-3$. Therefore, by condition $B_0$ we have,
$$
d(y_1)\geq 2a-2 \quad \hbox{and} \quad d(y_2)\geq 2a-2  \eqno (6) 
 $$
since $\{u_3,y_1,y_2\}\rightarrow u_4$ by (5).

Assume that $k\geq 3$, i.e., the cycle $C$ has length at least 6. Since $D$ contains the following arcs $x_1\rightarrow u_1, u_4\rightarrow y_1, u_4\rightarrow u_1$ and $y_1\rightarrow u_4$, it is not difficult to show that $y_{2}x_4\notin A(D)$ and $u_{2}y_2\notin A(D)$ 
   (possibly $x_4=x_1$) (for otherwise, if $y_{2}\rightarrow x_4$, then the cycle 
$x_1u_1u_2u_3u_4y_1x_2y_2x_4\ldots x_1$ has length $2a-2$; if $u_{2}\rightarrow y_2$, then now the cycle
$x_1y_1u_4u_1u_2y_2\ldots  x_1$ has length $2a-2$, which is a contradiction).
This together with  $u_{4} y_2\notin A(D)$ (by (2)) gives $d(y_2)\leq 2a-3$. This  contradicts the second inequality of (6).

Now assume that $k=2$. It is easy to see that  $y_{1}x_1\notin A(D)$ (for otherwise, $y_1x_1u_1u_2u_3u_4y_1$ is a cycle of length 6). This,  $u_{2}y_1\notin A(D)$ and $d(y_1)\geq 2a-2$ (by (6)) imply that $y_1\rightarrow u_2$ and $x_2\rightarrow y_1$. 
For the vertex $y_2$ we have that $u_4y_2\notin A(D)$, $y_2 u_2\notin A(D)$ and $x_1\rightarrow y_2\rightarrow x_2$, $u_2\rightarrow y_2$ because of the second inequality of (6). Therefore, $y_2x_1y_1u_4u_1u_2y_2$ is a cycle of length 6, a contradiction.

\textbf{Subcase 2.2}. $d(x_1)\geq 2a-2$.

By (2) we have $x_{1} u_3$ and $y_{k-1} x_1\notin A(D)$. This together with $d(x_1)\geq  2a-2$ implies that  $d(x_1)= 2a-2$ and 
$$
u_1\rightarrow x_1,\, \, u_3\mapsto x_1,\, \,  x_1\rightarrow y_k . \eqno (7)
$$

Assume that $k\geq 3$. Then, because of $y_{k-1} x_1\notin A(D)$, $y_{k-1}\not=y_1$ and (7), we have 
$$
\{u_1,u_3,y_1\}\rightarrow x_1.  
$$

It is not difficult to show that $d^+(y_k,\{u_2,u_4\})=0$ (since $\{u_1,u_3\}\rightarrow x_1$ and $u_4 \rightarrow u_1$),
 and 
$x_{k-1} y_k\notin A(D)$ (if $x_{k-1}\rightarrow y_k$, then $x_1u_1u_2u_3u_4y_1\ldots x_{k-1}y_kx_1$ is a cycle of length $2a-2$). 
Therefore, $d(y_k)\leq 2a-3$. Then from  condition $B_0$ we have $d(u_1)\geq 2a-2$ and $d(y_1)\geq 2a-2$, 
since
$\{u_1,y_1,y_k\}\rightarrow x_1$.  The inequality $d(u_1)\geq 2a-2$ together with $u_{1} u_4$ and $x_{k} u_1\notin A(D)$ implies that $u_1\leftrightarrow x_2$. Now  we have that the arcs $y_{1}u_4$, $u_{2} y_1$ and $y_{1} x_3$ are not in $A(D)$. Thus $d(y_1)\leq 2a-3$, which contradicts that $d(y_1)\geq 2a-2$.

Now assume that $k=2$. We consider the vertex $y_2$. It is easy to see that $y_{2} u_2\notin A(D)$ and $d(y_2,\{u_4\})=0$ (if $y_{2}\rightarrow u_4$, then, by (7), we have that $y_2u_4u_1x_1y_1x_2y_2$) is a cycle of length 6. Therefore, $d(y_2)\leq 2a-3$. This together with  
 condition $B_0$ implies that $d(u_3)\geq 2a-2$, since $\{y_2,u_3\}\rightarrow x_1$. The inequality $d(u_3)\geq 2a-2$ together with $x_{1} u_3$, $u_{3} x_2\notin A(D)$ implies that $x_{2}\rightarrow u_3$. Hence, the cycle $x_2u_3u_4u_1x_1y_1x_2$ is a cycle of length 6, a contradiction. Lemma 3.4 is proved. \fbox \\\\

\section {The proofs of the main results}

\textbf{Proof of Theorem 1.8.} Suppose, on the contrary, that $D$ is not a directed cycle and $D$ contains no cycle of length $2a-2$. Let $C=x_1y_1x_2y_2\ldots x_ky_kx_1$ be a longest non-Hamiltonian cycle in $D$, where $x_i\in X$ and $y_i\in Y$. By Lemma 3.2, $D$ contains a non-Hamiltonian cycle of length at least 4, i.e., $2\leq k\leq a-2$.
 By Lemma 3.1(ii), $D$ contains a $C$-bypass. Let $P:=uu_1\ldots u_sv$ be a $C$-bypass ($s\geq 1$). Suppose that the gap of $P$ is minimum  among the gaps of all $C$-bypasses.\\
 
From Lemma 3.4 it follows that $|V(C[u,v])|\geq 3$. The proof splits into cases, depending on the length of $C$-bypass $P$
and on the length of $C[u,v]$.

\textbf{Case 1.} $s\geq 2$.

Then $2\leq s\leq  |V(C[u,v])|-2$. Note that $\{u_s,v^-_C\}$ is a dominating pair. 
Since $C$ is a longest non-Hamiltonian cycle in $D$ and $C$-bypass $P$ has the minimum gap among the gaps of all $C$-bypasses, it follows that $v^-_C$ is not adjacent to any vertex on 
$P[u_1,u_s]$ and $u_s$ is not adjacent to any vertex on $C[u^+_C,v^-_C]$. Notice that each of 
$P[u_1,u_s]$ and $C[u^+_C,v^-_C]$ contains at least one vertex from each partite set. Therefore, 
$max\{d(u_s),d(v^-_C)\}\leq 2a-2$, which contradicts condition $B_1$ since $\{u_s,v^-_C\}$ is a dominating pair.\\ 

\textbf{Case 2.} $s= 1$.

It is easy to see that $u$ and $v$ belong to the same partite  set and $|V(C[u,v])|$ is odd. Now, without loss of generality, assume  that $x_1=u$,  
$v=x_r$ and $u_1=y$. Then $y\in Y$. Denote $R:=V(D)\setminus V(C)$ and $C':=V(C[y_1,y_{r-1}])$. Now we consider two subcases
($r\geq 3$ and $r=2$).

\textbf{Subcase 2.1.} $r\geq 3$.

Let $x$ be an arbitrary vertex of $ X\cap R$. 
Note that $\{y,y_{r-1}\}$ is a dominating pair. Recall that $y$ is not adjacent to any vertices of $C'$. Therefore,
$$
d(y)\leq 2a-2 \quad  \hbox{and} \quad  d(x_i)\leq 2a-2 \quad \hbox{for all} \quad x_i\in X\cap C',  \eqno (8)
$$
since $X\cap C'\not= \emptyset$. The first inequality of (8) together with  condition $B_1$ implies that
$$
d(y_{r-1})\geq 2a-1,  
$$
which in turn  implies that $y_{r-1}$ and every vertex of $X$ are adjacent. In particular, $y_{r-1}$ and $x$ are adjacent.

Assume first that $x\rightarrow y_{r-1}$. Then $\{x,x_{r-1}\}$ is a dominating pair. By the second inequality of (8) we have 
$d(x_{r-1})\leq 2a-2$. 
On the other hand, since $C$-bypass $P$ has the minimum gap among the gaps of all $C$-bypasses, it follows that $y_1 x$ and $y x\notin A(D)$. Therefore, $d(x)\leq 2a-2$. Thus,
$max\{d(x),d(x_{r-1})\}\leq 2a-2$,
 which is a contradiction since $\{x,x_{r-1}\}$ is a dominating pair.\\ 

Assume second that  $x y_{r-1}\notin A(D)$. Then $y_{r-1}\rightarrow x$, since $x$ and $y_{r-1}$ are adjacent. By the arbitrariness of $x$, we may assume that $y_{r-1}\mapsto X\cap R$. Combining this with $d(y_{r-1})\geq 2a-1$ we obtain that $|R|=2$, i.e., $|V(C)|=2a-2$, which contradicts the supposition that $D$ contains no cycle of length $2a-2$.\\

\textbf{Subcase 2.2.} $r=2$.

Now note that $\{y,y_1\}$ is a dominating pair. By condition $B_1$ ,
$max\{d(y),d(y_{1})\}\geq 2a-1$.
Without loss of generality, we may assume that $d(y)\geq 2a-1$ (if $d(y_1)\geq 2a-1$, then we will consider the cycle
 $Q:=x_1yx_2\ldots x_1$, which has the same length as $C$ and $x_1y_1x_2$ is a $Q$-bypass), which in turn implies that $y$ and every vertex of $X$ are adjacent. 
Notice that $|X\cap R|\geq 2$, since $|R|\geq 4$.
This together with $d(y)\geq 2a-1$ implies that there exists a vertex $x\in X\cap R$ such that $x\leftrightarrow y$. To complete the proof we now will prove the following three claims below.\\ 

\textbf{Claim 1.} For any $x_i\in X\cap V(C)$, if $x_i\rightarrow y$, then $x y_i\notin A(D)$; if 
$y\rightarrow x_i$, then $y_{i-1} x\notin A(D)$.

Indeed, for otherwise if $x_i\rightarrow y$ and  $x\rightarrow y_i$, then the cycle $x_iyxy_ix_{i+1}\ldots x_1$ is longer than $C$, and if 
$y\rightarrow x_i$ and $y_{i-1}\rightarrow x$, the cycle $y_{i-1}xyx_iy_i\ldots y_{i-1}$ is longer than $C$, a contradiction to the maximality of $C$. \fbox \\\\

From Claim 1 it immediately follows the following :

\textbf{Claim 2.} If there exists a vertex $x_i\in X\cap V(C)$ such that $x_i\rightarrow y\rightarrow x_{i+1}$, then $x$ and $y_i$ are not adjacent. \fbox \\\\

\textbf{Claim 3.} If there is a vertex 
$x_i\in X\cap V(C)$ such that $x_i\leftrightarrow y$, then 
$x_{i+1}\mapsto y$ and $y\mapsto x_{i-1}$ are impossible.

\textbf{Proof of Claim 3.} Suppose, on the contrary, that  there is a vertex $x_i\in X\cap V(C)$ such that $x_i\leftrightarrow y$ and 
 $x_{i+1}\mapsto y$ or  $y\mapsto x_{i-1}$. 
Combining this  with $d(y)\geq 2a-1$, we obtain that $y$ and every vertex of $X\cap R$ form a 2-cycle. In particular, any two vertices of $X\cap R$ form a dominating pair.
 If $x_{i+1}\mapsto y$, 
then $y\rightarrow x_{i+2}$ and, by Claim 2, the vertex $y_{i+1}$ and every vertex $z$ of $X\cap R$  are nonadjacent. If $y\mapsto x_{i-1}$, then $y\rightarrow x_{i+1}$ and, by Claim 2, the vertex $y_{i}$ and every vertex $z$ of $X\cap R$ are 
nonadjacent. In both cases we have $d(z)\leq 2a-2$ for all $z\in X\cap R$, which contradicts  condition $B_1$ since
 $|X\cap R|\geq 2$. Claim 3 is proved. \fbox \\\\

 Now we can finish the proof of Theorem 1.8.

From Claim 3 and $d(y)\geq 2a-1$ it follows that $y$ and every vertex of $X\cap V(C)$ form a 2-cycle. Combining this with Claim 2, we obtain that $x$ and every vertex $y_i$ are not adjacent, which in turn implies that 
$$
d(y_i)\leq 2a-2 \quad \hbox{and} \quad d(x)\leq 2a-2. \eqno (9)
$$
Note that for any $x_i$, $\{x,x_i\}$ is a dominating pair, since $\{x,x_i\}\rightarrow y$. This together with the second inequality of (9) and condition $B_1$ implies that  $d(x_i)\geq 2a-1$ for all $x_i$. 
Hence, any $x_i$ forms a 2-cycle with every vertex of $Y$ maybe except one. If $x_i$ and $y_j$ form a 2-cycle and $y_j\not=y_{i-1}$,then $\{y_{i-1},y_j\}\rightarrow x_i$ 
(i.e., $\{y_{i-1},y_j\}$ is a dominating pair) which is a contradiction because of (9). 
Assume therefore that $d(x_i,\{y_j\})\leq 1$ for all $y_j\not= y_{i-1}$. 
It follows that $y_i x_i\notin A(D)$ and $|Y\cap V(C)|=2$, i.e., the cycle $C$ has length to equal 4.
Thus, we may assume that for any $x_i$ there is no $y_j$ other than $y_{i-1}$ such that $y_j\rightarrow x_i$. 
We have $y_1x_1\notin A(D)$, $y_2 x_2\notin A(D)$, $x_2\rightarrow y_1$ and $x_1\rightarrow y_2$. 
Therefore $x_2$ forms a 2-cycle with every vertex of $Y\cap R$ since $d(x_2)\geq 2a-1$ and $y_2x_2\notin A(D)$.
Let $w\in Y\cap R$ be an arbitrary vertex other that $y$.
 Then $w\leftrightarrow x_2$ and hence, $\{y_{1},w\}$ is a dominating pair. By the first inequality of (9) and  condition $B_1$ we have $d(w)\geq 2a-1$, which in turn implies that 
$x\rightarrow w$ or $w\rightarrow x$. If $x\rightarrow w$, then $x_1yxwx_2y_2x_1$ is a cycle of length 6. If $w\rightarrow x$, then $x_2wxyx_1y_1x_2$ is a cycle of length 6. In both cases we have a longer non-Hamiltonian cycle than $C$, which is a contradiction and completes the discussion of Subcase 2.2. Theorem 1.8 is proved. \fbox \\\\

\textbf{Proof of Theorem 1.9.} Let $C=x_1y_1x_2y_2\ldots x_{a-1}y_{a-1}x_1$ be a cycle of length $2a-2$ in $D$, where $x_i\in X$ and $y_i\in Y$. Let the vertices $x$ and $y$ are not on $C$, where $x\in X$ and $y\in Y$. We distinguish two cases.

\textbf{Case 1.} $d(x)\geq 2a-2$ or $d(y)\geq 2a-2$.

Without loss of generality, we may assume that $d(x)\geq 2a-2$.  Since $d(x)=d(x,\{y\})+d(x, V(C))$ and since $a\geq 4$, it follows that 
$d(x, V(C))\geq 2a-4\geq a$. Therefore, by Lemma 3.3, $D$ contains cycles of all even lengths less than or equal to $2a-2$.\\

\textbf{Case 2.} $d(x)\leq 2a-3$ and   $d(y)\leq 2a-3$.

Without loss of generality, we may assume that $x\rightarrow y_1$ since $D$ is strong. We have, $\{x,x_1\}\rightarrow y_1$, i.e.,  $\{x,x_1\}$ is a dominating pair. Therefore,  by condition $B_0$, $d(x_1)\geq 2a-2$. It is easy to see that $D$ contains a cycle of length two since $a\geq 4$.

Assume first that there exists a vertex $y_l$ of $V(C)$ such that $x_1$ and $y_l$ are not adjacent. It is clear that $y_l\notin 
\{y_1,y_{a-1}\}$. Moreover, $x_1$ together with every vertex $y_i$ other than $y_l$ forms a 2-cycle. Therefore $x_1y_1x_2\ldots y_jx_1$ is a cycle of length $2j$ for every $j\in [1,l-1]\cup [l+1,a-1]$. 

Now we will show that $D$ contains also a cycle of length $2l$. If $3\leq l\leq a-3$, then $x_1y_2x_3\ldots x_{l+1}y_{l+1}x_1$ is a cycle of length $2l$.
Assume therefore that $l=a-2$ or $l=2$. Let $l=a-2\geq 3$. Then $x_1\leftrightarrow y_2$ and $x_1y_2x_3\ldots y_{a-1}x_1$
is a cycle of length $2l=2a-4$. If $l=a-2=2$, then $a=4$ and $x_1\leftrightarrow y$. 
By condition $B_0$, it is not difficult to see that 
$d(y_3)\geq 2a-2$ and $d(y_1)\geq 2a-2$ since $\{y,y_1,y_3\}\rightarrow x_1$ and $d(y)\leq 2a-3$. From this it follows that $d(y_3,\{x,x_2\})\geq 2$, and in any possible case $D$ contains a cycle of length 4. If $l=2\not=a-2$, then $x_1\leftrightarrow y_{a-2}$ and $x_1y_{a-2}x_{a-1}y_{a-1}x_1$ is a cycle of length 4. Thus, $D$ contains cycles of every length
$2,  4, \ldots , 2a-2$.

Assume second that $x_1$ and every vertex $y_i$ of $V(C)$ are adjacent. If
$$
x_1\rightarrow \{y_2,y_3,\ldots , y_{a-2}\} \quad \hbox{or} \quad \{y_2,y_3,\ldots , y_{a-2}\}\rightarrow x_1 \eqno (10)
$$
then it is easy to check that $D$ contains a cycle of length $2k$ for all $k\in [1,a-2]$. Assume therefore that (10) is not true. 
Using the fact that $d(x_1)\geq 2a-2$, it is not difficult to see that $y_1\rightarrow x_1\rightarrow y_{a-1}$.
Indeed, if $y_1 x_1\notin A(D)$, then from 
$$
2a-4\leq d(x_1,V(C))=d^+(x_1,\{y_1,y_2,\ldots ,y_{a-1}\})+d^-(x_1,\{y_2,y_3,\ldots ,y_{a-1}\})
$$
it follows that 
$d^+(x_1,\{y_1,y_2,\ldots ,y_{a-1}\})=a-1$ or $d^-(x_1,\{y_2,y_3,\ldots ,y_{a-1}\})=a-2$,  which contradicts the assumption that (10) is not true. Therefore, $y_1\rightarrow x_1$. Similarly one can show that $x_1\rightarrow y_{a-1}$. Since $x_1$ and every vertex $y_i\in V(C)$ are adjacent and since (10) is false, using $y_1\rightarrow x_1\rightarrow y_{a-1}$, we obtain that
 there are two distinct vertices $y_k$ and $y_l$ such that $y_k\mapsto x_1$ and $x_1\mapsto y_l$ where $2\leq k\leq a-2$ and  $2\leq l\leq a-2$.
This together with $d(x_1)\geq 2a-2$ implies that $x_1$ and every vertex of $Y\setminus \{y_k,y_l\}$ form a 2-cycle. 

Assume first that $l<k$. Then $2\leq l<k\leq a-2$ and 
 $$
 \{y_1,y_2,\ldots , y_{l-1},y_{l+1},\ldots , y_{a-1}\}\rightarrow x_1 \quad \hbox{and} \quad 
x_1\rightarrow \{y_1,y_2,\ldots , y_{k-1}\}
$$
Therefore, the cycle $x_1y_1x_2\ldots x_jy_jx_1$, where $j\in [1,a-1]\setminus \{l\}$, is a cycle of length $2j$, and the cycle 
$x_1y_2x_3\ldots x_{l+1}y_{l+1}x_1$ is a cycle of length $2l$. Thus $D$ contains cycles of every length $2, 4, \ldots , 2a-2$.

 Assume next that $k<l$. Then $2\leq k<l\leq a-2$ and
 $\{y_1,y_2,\ldots , y_{l-1},y_{l+1},\ldots , y_{a-1}\}\rightarrow x_1 $.
From this we have that $x_1y_1x_2\ldots x_{j}y_{j}x_1$ is a cycle of length $2j$ for all $j\in [1,l-1]\cup [l+1,a-2]$. 

Now we want to show that $D$ contains also a cycle of length $2l$. If $k\geq 3$, then $x_1\rightarrow y_2$ and the cycle $x_1y_2x_3\ldots x_{l+1}y_{l+1}x_1$ is a cycle of length $2l$. Assume therefore that $k=2$. Then
$$
x_1\rightarrow \{y_3,y_4,\ldots , y_{a-1}\} \quad \hbox{and} \quad \{y_{l+1},y_{l+2},\ldots , y_{a-1}\}\rightarrow x_1.
$$
Hence, if $l\leq a-3$, then the cycle $x_1y_3x_4\ldots x_{l+2}y_{l+2}x_1$ is a cycle of length $2l$. Thus we have $k=2$  and $l=a-2$. It remains to show that $D$ contains a cycle of length $2l=2a-4$.
Recall that $y_2\mapsto x_1$ and $x_1\mapsto y_{a-2}$. Since $x_1\leftrightarrow y$, it follows that $\{y,y_{a-1}\}$ is a dominating pair. By the assumption of Case 2, $d(y)\leq 2a-3$. 
Therefore, from condition $B_0$ it follows that $d(y_{a-1})\geq 2a-2$. 
If $y_{a-1}\rightarrow x_2$ or $x_{a-2}\rightarrow y_{a-1}$, then $y_{a-1}x_2y_2\ldots x_{a-1}y_{a-1}$ or $x_{a-2}y_{a-1}x_1y_1\ldots y_{a-3}x_{a-2}$ is a cycle of length $2a-4$, respectively for $y_{a-1}\rightarrow x_2$ and $x_{a-2}\rightarrow y_{a-1}$. 
Assume therefore that 
$y_{a-1}x_2\notin A(D)$ and $x_{a-2} y_{a-1}\notin A(D)$. Then, since $d(y_{a-1})\geq 2a-2$, we have  $y_{a-1}\leftrightarrow x$ and hence, the cycle $y_{a-1}xy_1x_1y_3\ldots x_{a-1}y_{a-1}$ is a cycle of length $2a-4$. Thus we have shown that $D$ contains cycles of every even length $m$, $2\leq m\leq 2a-2$. Theorem 1.9 is proved. \fbox \\\\

\textbf{Proof of Theorem 1.10.} By Theorem 1.7, either $D$ contains a cycle of length $2a$ (i.e., $D$ is Hamiltonian) or $D$ is isomorphic  to the digraph $D(8)$. By Theorem 1.8, if $D$ is not a directed cycle, then it contains a cycle of length $2a-2$. Now by Theorem 1.9, $D$ contains cycles of every even length $2k$, where $1\leq k\leq a-1$. Therefore, if $D$ is other than $D(8)$ and a directed cycle, then $D$ is even pancyclic. Theorem 1.10 is proved.\\

\end{document}